\documentclass[11pt]{article}
\usepackage{amssymb}
\usepackage{}
\usepackage{pifont}
\usepackage{amscd}
\usepackage{amsfonts}
\usepackage{amsmath, amsthm}
\usepackage{amssymb, amsxtra}

\usepackage{graphicx}

\usepackage{multicol}
\usepackage{epstopdf}
\usepackage{epsf}
\usepackage{float}
\usepackage{extarrows}
\usepackage[a4paper, text={.8\paperwidth,.83\paperheight}, ratio=1:1]{geometry}
\usepackage[format=hang,font=small,textfont=it]{caption}
\usepackage{authblk}

\usepackage[colorlinks=true, linkcolor=blue, citecolor=blue]{hyperref}
\usepackage[english]{babel}
\usepackage{latexsym}

\usepackage{tikz}
\usetikzlibrary{calc}
\usepackage[tight]{subfigure}

\usepackage{color}
\definecolor{greenbean}{RGB}{199,237,204}

\newtheorem{thm}{Theorem}[]

\newtheorem{lemma}[thm]{Lemma}

\newtheorem{eg}{Example}[]

\newtheorem{Def}[thm]{Definition}

\def \ta{\tau}

\def \ta1{\tau_1}

\def \prodl{\prod\limits}

\begin{document}
\title{On Galois covers of a union of Zappatic surfaces of type $R_k$
\footnotetext{\hspace{-1.8em} Email address: Meirav Amram: meiravt@sce.ac.il; Cheng Gong: cgong@suda.edu.cn; 
Jia-Li Mo:  mojiali0722@126.com;\\
2020 Mathematics Subject Classification. 05E15, 14J10, 14J25, 14N20.}}

\author[1]{Meirav Amram}
\author[2]{Cheng Gong}
\author[2]{Jia-Li Mo}

 \affil[1]{\small{Shamoon College of Engineering, Ashdod, Israel}}
\affil[2]{\small{Department of Mathematics, Soochow University, Shizi RD 1, Suzhou 215006, Jiangsu, China}}

\date{}
\maketitle

\abstract{We investigate the topological structures of Galois covers of a union of two Zappatic surfaces of type $R_k$. We prove that the Galois covers of such surfaces are simply connected surfaces of general type. We also compute their Chern numbers and topological indices.}

\section{Introduction}\label{outline}


The moduli space of algebraic varieties is one of the interesting objects in geometry and topology. The moduli space for curves was constructed in the 1960s. In \cite{Gie}, the author constructed the moduli space for surfaces of general type over $\mathbb{C}$.
Then, classifying algebraic surfaces and studying their moduli space became one of the most thoroughly investigated subjects in algebraic geometry and topology; for example, see \cite{C1, C2, Hu1, Hu2}. One way to a classification of algebraic surfaces is to study generic projections. Take surface $X$ in $\mathbb{CP}^n$ and project it generically onto a plane. Now, let $S$ be the branch curve of such projection; it is almost completely determined to recover surface $X$ from the fundamental group $\pi_1(\mathbb{CP}^2-S)$ (see 
 \cite{KT}), but computations are complicated.
It is well-known that the Galois cover 
$X_{Gal}$ attached to a generic projection from $X$ is again a smooth projective surface. Its fundamental group $\pi_1(X_{Gal})$ is a quotient of $\pi_1(\mathbb{CP}^2-S)$, and therefore it is
expected to contain non-trivial information about the surface $X$. 
Bogomolov conjectured that an algebraic surface of general type with a positive index has an infinite fundamental group. There is a counter-example for this conjecture in \cite{MoTe87}, but the interesting thing is to see what happens when the topological index is negative.

Now, we introduce the objects that are studied in this paper and give important references to them. In \cite{zg2, Zappa}, Zappa first studied a type of surface degeneration, now known as Zappatic degeneration. Calabri, Ciliberto, Flamini, and Miranda classified the singularity of Zappatic degeneration in \cite{C-C-F-M-5, C-C-F-M-4}. Three types of Zappatic singularities are  $E_k$, $R_k$, and $S_k$, and they are very well studied and explained in those works.

In \cite{AGM}, we study Zappatic degenerations with only one singularity of type $R_k$ (for every $k \geq 3$) and classify their Galois covers.
In this paper, we glue two Zappatic surfaces, each of which degenerates to one singularity of type $R_k$, and the resulting surface is called $X_0$ (see Section \ref{method} for more details). It can be degenerated into a union of the cones over the same rational curves. After computing  their Galois covers, we get the  main theorem as follows:
\begin{thm}\label {1}
If $X$  is a smooth algebraic surface in $\mathbb{CP}^n$, it can be degenerated into a union of the cones over the same rational curves. Then the Galois cover $X_{Gal}$ is a simply connected surface of a general type.
\end{thm}

The paper is organized as follows:  In Section \ref{method}, we give the method and the terminology related to the paper. Moreover, we show the construction of the degeneration for an algebraic surface. In Section \ref{Results}, we show the fundamental groups of the Galois covers of these surfaces and determine the results. In Section \ref{Chern}, we compute the Chern numbers of such surfaces and, after creating the topological indices, we present \begin{thm}
  The topological index of the above Galois cover in Theorem \ref{1} is negative.
\end{thm}

\paragraph{Acknowledgements:}
We thank Prof. Ciro Ciliberto for his useful discussions about the degeneration of surfaces. The work is partly supported by NSFC~(12331001)  and the Natural Science Foundation of Jiangsu Province~(BK 20181427). The work was written during the third co-author's visit to SCE.

\section{Method and terminology}\label{method}

In this section, we describe in detail the methods that we use to determine the fundamental group of the Galois cover of the surface that is a union of two Zappatic surfaces with one singularity of type $R_{k}, \ k \geq 4$. These are degenerated surfaces, corresponding to the original smooth projective algebraic surface, and are represented as $R_k \bigcup R_k$. The results for $k=1,2,3$ were proved in \cite{degree6,A-R-T}. The reader can refer to \cite{cog1} for the monodromy study, \cite{Libgober1, Paris} for the fundamental group study, \cite{Enri} for understanding the topology of our objects and arrangements, and \cite{C-C-F-M-5} for the study of degenerations.

The fundamental group of the Galois cover is a significant invariant of the surface, as explained in the introduction; we are going to calculate and determine that invariant. We explain now the methods to do it and recall information from \cite{AGM}, in which we computed the fundamental group of the Galois cover of the surface that degenerates to a cone over a smooth rational curve.

Let $ X $ be a projective algebraic surface embedded in projective space $ \mathbb{CP}^n $, for some $n$. A generic  projection  $f:\mathbb{CP}^n\to\mathbb{CP}^2$, restricted to $ f|_X $, is branched along a curve $ S\subset \mathbb{CP}^2 $. The branch curve $ S $ can tell a lot about $X$, but it is difficult to describe it explicitly. Therefore, a degeneration of $X$ can ease the work, especially because we can construct a sequence of partial degenerations $X:~=X_r \leadsto X_{r-1} \leadsto \cdots X_{r-i} \leadsto X_{r-(i+1)} \leadsto \cdots \leadsto X_0$. The degenerated surface $X_0$ is a union of planes, with each plane projectively equivalent to $\mathbb{CP}^2$. 

\begin{figure}[ht]
\begin{center}
\scalebox{0.5}{\includegraphics{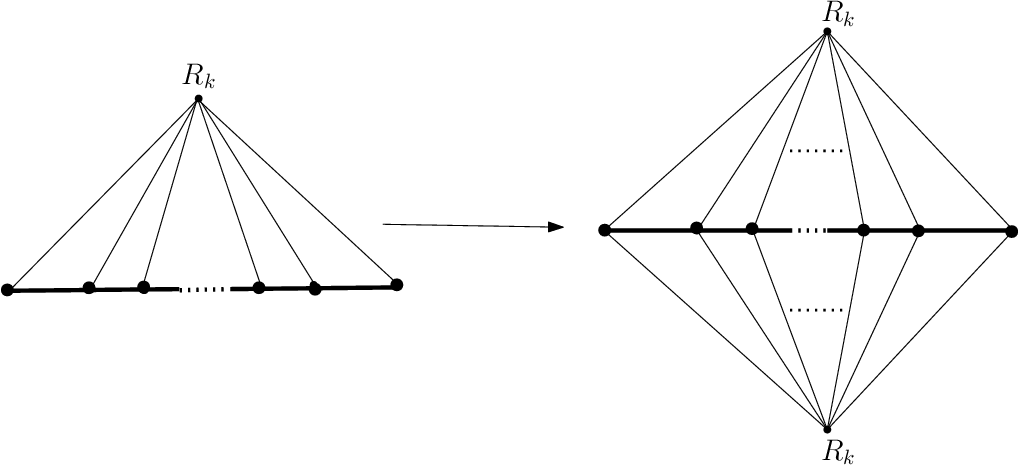}}
\end{center} \caption{Zappatic degenerations with Zappatic singularities of type $R_k$}\label{Rk-fig}
\end{figure}

On the left side of Figure \ref{Rk-fig} we have a degeneration with a unique Zappatic singularity of type $R_k, \ k \geq 3$.
The singular surface is the union of $k$ planes with a $k$-tuple point  (every triangle represents a plane).
The $k$-tuple point is a Zappatic singularity of type $R_k$. As such, the singular surface can be degenerated into a cone over a smooth rational curve. For the details see \cite{AGM}.

Then we can glue two such singular surfaces along the bold edges and get a singular surface  $X_0$ that has two Zappatic singularities of type $R_k$ on the right side of the figure. 
The basic degeneration of $X$ is a union of two rational normal scrolls of degree $k$ glued along a hyperplane section.
Next, we can use a method similar to Sec 2.2 in \cite{CLM} to construct a flat morphism: $f: X\to X_0$, where $X$ is a smooth surface. $X$ sits in $\mathbb{P}^{n+2}$, its degree is  $2k$, and the genus of the general hyperplane section is $k-1$.  However, $f$ is the degeneration we study in this work (for convenience of computations, in Subsection \ref{n+1-type} we write $k=n+1$).

\begin{figure}[ht]
\begin{center}
\scalebox{0.6}{\includegraphics{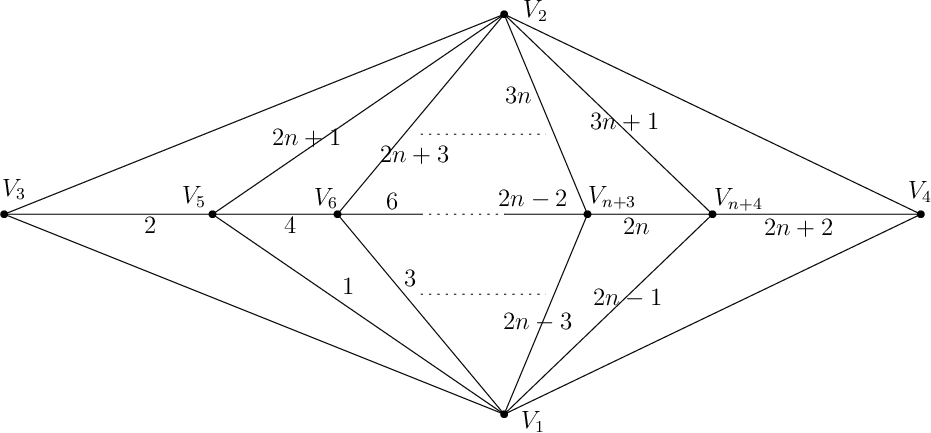}}
\end{center} \caption{Global ordering in the degeneration}\label{RnRn}
\end{figure}

Now we follow Figure \ref{RnRn}, and pay attention to vertices $V_1$ and $V_2$ that are both Zappatic singularities of type $R_{n+1}$ (each is an intersection of $n$ lines).

We number the other vertices $V_3,\dots,V_{n+4}$, as shown in Figure \ref{RnRn}. The lines in this degeneration are numbered as well, and we are going to have $3n+1$ of them. The collection of these $3n+1$ lines is the branch curve $S_0$ of $X_0$. 
In the figure are these two Zappatic singularities of type $R_{n+1}$ (vertices $V_1$ and $V_2$), two more vertices ($V_3$ and $V_4$) that have one line passing through each of them (line 2 passes through $V_3$ and line $2n+2$ passes through $V_4$), and the rest of the vertices (namely $V_5,\dots, V_{n+4}$) are intersections of four lines each. The local numbers of the four lines passing through each of the vertices $V_5, \dots, V_{n+4}$ are the same, and this similar numbering enables us, later on, to write generalized lemmas about the presentation of the fundamental group. 

We use a regeneration chain of branch curves $S_i$ ($0 \leq i \leq r$), of the generic projections $f_i: X_{i} \rightarrow \mathbb{CP}^2$ to regenerate $S_0$. The curve $S_r$ is the branch curve $ S $ and it is cuspidal with conics and lines, and with nodes and branch points as additional singularities. The degree of $S$ is double the one of $S_0$ because each line in $S_0$ regenerates somewhere along this regeneration chain to a conic, a pair, or parallel lines.

By the van-Kampen theorem for a cuspidal curve \cite{vk}, we can compute group $G:=\pi_1(\mathbb{CP}^2-S)$, which is the fundamental group of the complement of branch curve $S$ in $\mathbb{CP}^2$. 
To compute the presentation we first use the global numbering of the lines in the degeneration from Figure \ref{RnRn} and denote the standard generators of $G$ as $1, 1', \dots, {3n+1}, {(3n+1)'}$.
The relations arising from each vertex depend on the local ordering of the lines at this point, see Figure \ref{RnRn}.
Then we write a list of relations according to the types of singularities ($i, j$ are the standard generators in $G$): 
\begin{enumerate}
\item For every branch point of a conic, we have a relation of the form ${j}={j'}$, and
\item\label{vK2} for every node, we have $[i,j] = e$ where $i,j$ are the lines meeting in this node, and 
\item for every cusp, $\langle i,j \rangle=1$ where $i,j$ are two components involved in this cusp, we have
\item\label{vK4} the ``projective relation'' $\prodl_{i=(3n+1)}^1  {j'}\ j = e$, and
\item\label{par} a commutator relation for every pair of disjoint edges from Figure \ref{RnRn}, which meet in the branch curve.
\end{enumerate}

As a consequence, other than the projective relation, the relations can be computed locally. In other words, $G$ is an amalgamated product of ``local groups,'' one for each vertex, generated by the generators corresponding to the lines meeting in the vertex --- modulo the projective relation of type (\ref{vK4}) and the relations of type (\ref{par}). Our intermediate goal is to find group  $G_1:=G/\langle {j^2, j'^2}\rangle$, from which we define thereafter a surjection onto the symmetric group $S_n$ and approach our main goal --- the kernel of this map.

We start with the local relations that vertices $V_5,\dots,V_{n+4}$ contribute to group $G$.

\begin{lemma}\label{basic-abcd}
Figure \ref{inner4} depicts an intersection of four lines, being locally ordered as $ a < b < c < d$.

\begin{figure}[ht]
\begin{center}
\scalebox{0.8}{\includegraphics{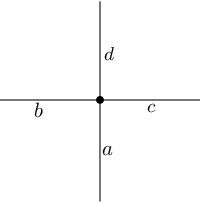}}
\end{center} \caption{Intersection of four lines}\label{inner4}
\end{figure}

A point that is an intersection of four lines contributes to $G$ the following list of relations:

\begin{equation}\label{eq1.1}
\langle a',b\rangle=\langle a',b'\rangle=\langle a',b^{-1}b'b\rangle=e,
\end{equation}
\begin{equation}\label{eq1.2}
\langle c,d\rangle=\langle c',d\rangle=\langle c^{-1}c'c,d\rangle=e,
\end{equation}
\begin{equation}\label{eq1.3}
[b'ba'b^{-1}b'^{-1},d]=e,
\end{equation}
\begin{equation}\label{eq1.4}
[b'ba'b^{-1}b'^{-1},c^{-1}c'^{-1}d^{-1}d'dc'c]=e,
\end{equation}
\begin{equation}\label{eq1.5}
\langle a,b\rangle=\langle a,b'\rangle=\langle a,b^{-1}b'b\rangle=e,
\end{equation}
\begin{equation}\label{eq1.6}
\langle c,d^{-1}d'd\rangle=\langle c',d^{-1}d'd\rangle=\langle c^{-1}c'c,d^{-1}d'd\rangle=e,
\end{equation}
\begin{equation}\label{eq1.7}
[b'bab^{-1}b'^{-1},d^{-1}d'd]=e,
\end{equation}
\begin{equation}\label{eq1.8}
[b'bab^{-1}b'^{-1},c^{-1}c'^{-1}d^{-1}d'^{-1}dd'dc'c]=e,
\end{equation}
\begin{equation}\label{eq1.9}
b'ba'ba'^{-1}b^{-1}b'^{-1}=dc'd^{-1},
\end{equation}
\begin{equation}\label{eq1.10}
b'ba'b'a'^{-1}b^{-1}b'^{-1}=dc'cc'^{-1}d^{-1},
\end{equation}
\begin{equation}\label{eq1.11}
b'baba^{-1}b^{-1}b'^{-1}=d^{-1}d'dc'd^{-1}d'^{-1}d,
\end{equation}
\begin{equation}\label{eq1.12}
b'bab'a^{-1}b^{-1}b'^{-1}=d^{-1}d'dc'cc'^{-1}d^{-1}d'^{-1}d.
\end{equation}
\end{lemma}

In addition, vertices $V_3$ and $V_4$ contribute the relations $2=2'$ and $(2n+2)=(2n+2)'$ to $G$, because lines $2$ and $2n+2$ regenerate to conics and the van-Kampen theorem supplies such relations.

We have also the relations coming from the two Zappatic singularities of type $R_{n+1}$, these are vertices $V_1$ and $V_2$. Each vertex contributes local relations, as listed in Lemma \ref{Rn-lemma1}, below. The relations in Lemma \ref{Rn-lemma1} are those simplified rejections from Subsection \ref{claim}, with the condition $j^2 = e, \ j'^2 = e$ included as well. This means that the relations are the ones that will appear later on in $G_1$. 

For the convenience of readers, we refer to the relevant material from \cite{AGM} without proof, thus making our exposition self-contained. In \cite{AGM}, we focus on the Zappatic degenerations with only one Zappatic singularity (e.g.,~$O$ of type $R_{n+1}$ for $n \geq 2$), and compute the fundamental group $G=\pi_1 (\mathbb{CP}^2-S)$ of the complement of the branch curve of a Zappatic surface.

We provide full relations given by all singularities in the fundamental group $G$. As to our prerequisites for understanding this paper, we expect our audience to be familiar with the simplest relations obtained by one Zappatic singularity. To find some rules directly, we provide a direct demonstration for relations and a short proof.

\begin{figure}[ht]
\begin{center}
\scalebox{0.8}{\includegraphics{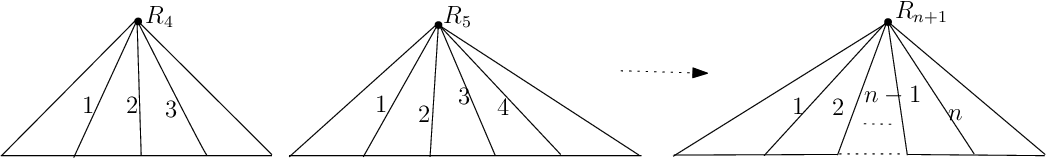}}
\end{center} \caption{A series of Zappatic singularity of type ${R_{n+1}}$}\label{R4R5Rn}
\end{figure}

Take a series of Zappatic singularities of type $R_k, \ k=4,5,\cdots,n+1$, see Figure \ref{R4R5Rn}. 

\begin{lemma}\label{Rn-lemma1}
Given an $R_{n+1}$ Zappatic singularity with lines $1,2,\dots,n$, as in Figure \ref{R4R5Rn}-third.

The relations in $G_1$ are: 
\begin{eqnarray}\label{eq1.19}
&~\langle 1',2\rangle=\langle 1',2'\rangle=\langle 1',22'2\rangle=e;~\langle 2',3\rangle=\langle 3,2'22'1'2'22'\rangle=\langle 2,3\rangle=e,\\
&\langle 3,4\rangle=\langle3',4\rangle=e; \cdots; ~~\langle n-2,n-1\rangle=\langle(n-2)',n-1\rangle=e,\\
&\langle n-1,n\rangle=\langle(n-1)',n\rangle=e,\\
&~1=2'21'22';~~~33'3=2'2322';~~~44'4=3'34'33'; ~~\cdots,\\
&(n-1)(n-1)'(n-1)=(n-2)'(n-2)(n-1)'(n-2)(n-2)', \\
& nn'n=(n-1)'(n-1)n'(n-1)(n-1)',\\
&~[1,3]=[1,3']=[1',3]=[1',3']=e,\\
&~[1,4]=[1',4]=e;~~ \cdots;~~ [1,n-1]=[1',n-1]=e; ~~[1,n]=[1',n]=e,\\
&[2,4]=[2',4]=e; ~~\cdots; ~~[2,n-1]=[2',n-1]=e; ~~[2,n]=[2',n]=e,\\
&[3,5]=[3',5]=e; ~~\cdots; ~~[3,n-1]=[3',n-1]=e; ~~[3,n]=[3',n]=e,\\
&\cdots\cdots\cdots\\
&[n-2, n]=[(n-2)', n]=e.
\end{eqnarray}
\end{lemma}

\begin{proof}
Zappatic singularity of type $R_4$ contributes the following relation in $G_1$:

\begin{equation}\label{5.1}
\langle 1',2\rangle=\langle 1',2'\rangle=\langle 1',22'2\rangle=e,
\end{equation}
\begin{equation}\label{5.2}
1=2'21'22',
\end{equation}
\begin{equation}\label{5.3}
\langle 2'21'21'22',3\rangle=\langle 2'21'2'1'22',3\rangle=\langle 2'21'22'21'22',3\rangle=e,
\end{equation}
\begin{equation}\label{5.4}
3'=32'21'2'21'22'32'21'22'1'22'3,
\end{equation}
\begin{equation}\label{5.5}
[1,3]=[1,3']=[1',3]=[1',3']=e.
\end{equation}

These relations can be simplified as follows:

Because some of these relations are already the simplest, we return to relations (\ref{5.3}) and (\ref{5.4}). For these, substituting (\ref{5.1}), (\ref{5.2}), and (\ref{5.5}) into them, we obtain the simplest relations which are $\langle 2',3\rangle=\langle3,2'22'1'2'22'\rangle=\langle 2,3\rangle=e$ and $33'3=2'2322'$. To prove completeness, here we provide an example of the detailed process showing how to get relation $\langle 2',3\rangle=e$:

\begin{itemize}
  \item  Firstly: $1=2'21'22'\Longrightarrow 1'=22'12'2$,
  \item  Secondly: substituting $1'=22'12'2$ into $\langle 2'21'22'21'22',3\rangle=e\Longrightarrow\langle 12'1,3\rangle=e$,
  \item  Thirdly: using $[1,3]=e$, we get relation $\langle 2',3\rangle=e$.
\end{itemize}
Finally, we get the simplest relations in the following:
\begin{eqnarray}\label{eq2.17}
&~\langle 1',2\rangle=\langle 1',2'\rangle=\langle 1',22'2\rangle=e;~~~\langle 2',3\rangle=\langle 3,2'22'1'2'22'\rangle=\langle 2,3\rangle=e,\\
&~1=2'21'22';~~~33'3=2'2322';~~[1,3]=[1,3']=[1',3]=[1',3']=e.
\end{eqnarray}

Relations associated with $R_5$ are the following:
\begin{eqnarray}\label{eq1.17}
&~[1,434]=[1,43'4]=[1',434]=[1',43'4]=e,\\
&~\langle3,4\rangle=\langle3',4\rangle=\langle 33'3,4\rangle=e,\\
&~[2'2122',4]=[2'21'22',4]=e,\\
&~[2,4]=[2',4]=e,\\
&~[3'32'2122'33',44'4]=[3'32'21'22'33',44'4]=e,\\
&~[3'3233',44'4]=[3'32'33',44'4]=e,\\
&~44'4=3'34'33'.
\end{eqnarray}
And they will be simplified, in a similar way, to the following relations:
\begin{eqnarray}\label{eq1.18}
&~\langle 1',2\rangle=\langle 1',2'\rangle=\langle 1',22'2\rangle=e,\\
&~\langle 2',3\rangle=\langle 3,2'22'1'2'22'\rangle=\langle 2,3\rangle=e;~~\langle 3,4\rangle=\langle3',4\rangle=e,\\
&~1=2'21'22';~~~33'3=2'2322';~~~44'4=3'34'33',\\
&~[1,3]=[1,3']=[1',3]=[1',3']=e,\\
&~[1,4]=[1',4]=~[2,4]=[2',4]=e.
\end{eqnarray}

Now, looking at the simplified relations associated with $R_4$ and $R_5$ (and even to be convinced, in the case of $R_6$ and so on), we see the rule in the structure of the relations, which gives us the relations in the lemma. 
\end{proof}

\begin{figure}[ht]
\begin{center}
\scalebox{0.8}{\includegraphics{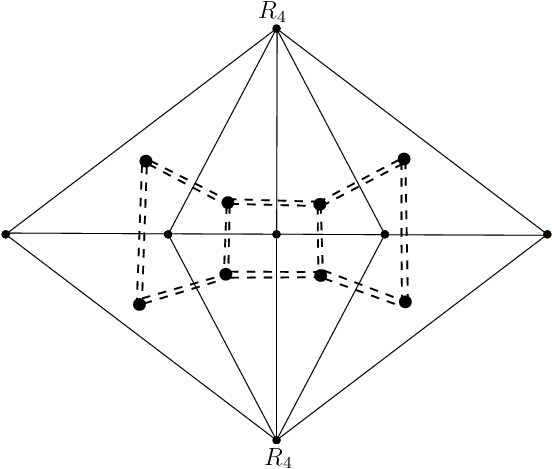}}
\end{center} \caption{The dual graph related to degenerated surface of type $R_4$ }\label{R4R4combine}
\end{figure}

The amalgamation is done by identifying the generators associated with the same line and combining all the above relations. This combining includes the projective relation and the commutations from type (\ref{par}) in the method.
When we collect all relations together, they include the dual graph related to $G_1$, which includes all pairs of generators  $j$ and $j'$. As an example, in Figure \ref{R4R4combine} we depict a union of two Zappatic surfaces with one singularity of type $R_{4}$. The double dashed lines everywhere represent those generators $j$ and $j'$ on the line $j$ in the degeneration. Later, all the above relations will be simplified further, under the condition $j^2 = e, \ j'^2 = e$, and we will get a simplified presentation of group $G_1$. By eliminating generators, we can get another dual graph, and ensure that the resulting simplified relations maintain the graph structure.
As an example of such a dual graph, we take the degeneration---the union of Zappatic surfaces with one singularity of type $R_4$ in Figure \ref{R4R4combine}. 
As it stands now, the group has $20$ generators that are $j, \ j'$ for $j=1, \dots, 10$. Thereafter, simplifications eliminate generators.

The algorithm of van-Kampen leads us eventually to compute the fundamental group of the Galois cover $X_{Gal}$ of $X$ as well.
\begin{Def}
We consider the fibered product arising from a general projection $f: X \to \mathbb{CP}^2$ of degree $n$ as
$$X\times_{f}\cdots\times_{f}X=\{(x_1, \ldots , x_n)\in X^n|~f(x_1)=\cdots=f(x_n)\}.$$
Let the extended diagonal be
$$ \triangle=\{(x_1, \ldots , x_n)\in X^n|~x_i=x_j, for ~some ~i\neq j\}.$$
The closure
$\overline{X\times_{f}\cdots\times_{f}X-\triangle}$ is called the Galois cover w.r.t. the symmetric group $S_n$ and denoted by $X_{Gal}$.
\end{Def}
We take a surjection from $G_1$ onto the symmetric group $S_n$, where each generator of $G_1$ maps to a transposition in the symmetric group $S_n$, according to the order of the lines and planes in the degeneration. We thus obtain a presentation of the fundamental group $\pi_1(X_{Gal})$ of the Galois cover, as the kernel of this surjection. Then we simplify the relations to produce a canonical presentation that identifies with $\pi_1(X_{Gal})$; it is hard difficult but the result is interesting and meaningful. 
In Section \ref{Results}, we give detailed simplifications and calculations of the group $\pi_1(X_{Gal})$.

\section{The main theorem}\label{Results}

We divide this section into three parts to elaborate our work. Firstly, for the convenience of calculations, we provide a lemma that will support the proof of the main theorem. 
Secondly, we provide two special cases of proof to prepare for the general results.
Finally, at the end of this section, we prove the main theorem for general $k$.

\subsection{Preliminary claims}\label{claim}

Figure \ref{inner4} depicts an intersection of four lines, being globally ordered as $ a < b < c < d$. A point that is an intersection of four lines contributes to $G$ the relations in Lemma \ref{basic-abcd}. Now, as we consider group $G_1$, it means that we add $j^2=e$ and $j'^2=e$.

\begin{lemma}\label{innerR4-lemma1}
If $b=b'$ (from Figure \ref{inner4}) is added to the relations from Lemma \ref{basic-abcd}, and we consider the conditions on $G_1$, then we have $c=c'$ and other relations as follows: 
\begin{equation}\label{4.1}
\langle a,b\rangle=e;~~\langle b,d\rangle=e;~~\langle d,c\rangle=e;~~\langle c,a\rangle=e,
\end{equation}
\begin{equation}\label{4.2}
[b,c]=e;~~[a,d]=e;~~a'=bdcdb;~~d'=dcabacd.
\end{equation}
\end{lemma}

\begin{proof}
When adding relation $b=b'$ to the relations (\ref{eq1.1})--(\ref{eq1.12}) in Lemma \ref{basic-abcd}, first, we have $c=c'$ by relations (\ref{eq1.9}) and (\ref{eq1.10}), and next, it leads to the following relations in $G_1$:
\begin{equation}\label{4.3}
\langle a',b\rangle=e;~~\langle c,d\rangle=e;~~[a',d]=e;~~[a',dd'd]=e, 
\end{equation}
\begin{equation}\label{4.4}
\langle a,b\rangle=e;~~\langle c',dd'd\rangle=e;~~[a,dd'd]=e;~~[a,dd'dd'd]=e,
\end{equation}
\begin{equation}\label{4.5}
a'=bdcdb;~~aba=dd'dcdd'd.
\end{equation}

Until now, it remains to prove relations $\langle b,d\rangle=e, \ \langle c,a\rangle=e, \ [b,c]=e, \ [a,d]=e, \ d'=dcabacd$. These missing relations can be obtained by relations (\ref{4.3})--(\ref{4.5}). The details are as follows: 
We conclude from $[a',dd'd]=e$ and $[a,dd'dd'd]=e$ that $[a,d'd]=e$ and $[a,d]=e$, hence, that $[a,d']=e$ by $[a,d'd]=e$ and $[a,d]=e$. Moreover, we get $[a',d']=e$ from $[a',d]=e$ and $[a',dd'd]=e$.  Substituting $a'=bdcdb$ into $\langle a',b\rangle=e$ and $[a',d]=e$, we obtain $\langle b,d\rangle=e$ and $[b,c]=e$ respectively. Then combining $\langle c',dd'd\rangle=e$ with $aba=dd'dcdd'd$, we have $d'=dcabacd$. Similarly, we substitute $a'=bdcdb$ and $d'=dcabacd$ into $[a',dd'd]=e, \ 
\langle c',dd'd\rangle=e$ and $[a,dd'dd'd]=e$, and obtain relations $[a,d]=e, \ \langle a,c\rangle=e$, and  $[c,d]=e$ respectively. 
\end{proof}


\subsection{A union of Zappatic surfaces with one singularity of type $R_4$($R_5$)} \label{6-type}

\begin{figure}[ht]
\begin{center}
\scalebox{0.8}{\includegraphics{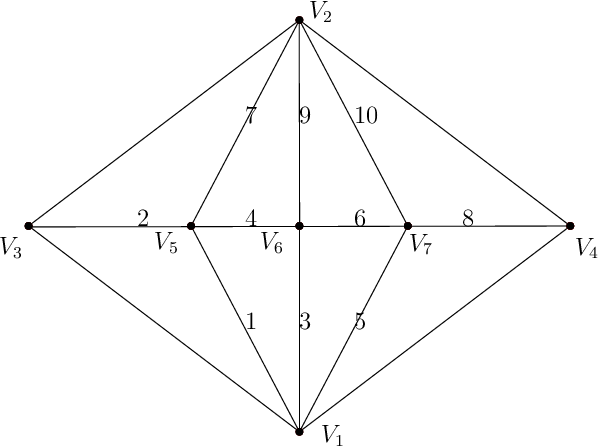}}
\end{center} \caption{The degenerated of the surface $R_4 \bigcup R_4$}\label{R4R4}
\end{figure}

Based on the above lemmas, we are in a position to compute the fundamental group of Galois covers of a union of Zappatic surfaces in this subsection. We first prove the basic form of the main theorem. In the following, we give a degeneration of a surface, which is a union of Zappatic surfaces, each of which with one singularity of type $R_4$ (notation for the surface will be $R_4 \bigcup R_4$) in Figure \ref{R4R4}. There are seven vertices $V_i (\ i=1,2,\dots,7$), and 10 lines which consist of the branch curve of the degenerated Zappatic surface.

\begin{thm}\label{thm44}
The Galois cover of a union of Zappatic surfaces, each of which with one singularity of type $R_4$, is simply connected.
\end{thm}

\begin{proof}
Explicitly, much information can be obtained directly from Figure \ref{R4R4}. Here $V_1$ and $V_2$ are Zappatic singularities of type $R_4$, with each vertex $V_i$ for $\ i=5,6,7$ is an intersection of four lines. Let $S_{R_4 \bigcup R_4}$ be the branch curve of $R_4 \bigcup R_4$.
The fundamental group $G=\pi_1(\mathbb{CP}^2-S_{R_4 \bigcup R_4})$ is generated by generators
${\{j,j'\}}^{10}_{j=1}$, as for now we work with $G_1=G/<j^2,j'^2>$.

Vertices $V_3$ and $V_4$ give rise to the relations $2=2'$ and $8=8'$ in $G_1$.
Having relation $2=2'$ in $V_5$, we get relation $4=4'$ by Lemma \ref{innerR4-lemma1}. Similarly, we get from $V_6$ that $6=6'$, and thereafter, we get from $V_7$ that $8=8'$.

Lastly, as we collect these resulting equalities, together with other relations, we note that relations related to $V_5, V_6$ and $V_7$ will be listed here (Lemma \ref{innerR4-lemma1}), and relations related to $V_1, V_2$ will be listed here as well (Lemma \ref{Rn-lemma1}). We have also the projective relation from (\ref{vK4}) and the commutations from (\ref{par}).

\begin{equation}\label{6.1}
2=2'; \ 4=4'; \ 6=6'; \ 8=8',
\end{equation}
\begin{equation}\label{6.2}
\langle 1,2\rangle=\langle 2,7\rangle=\langle 7,4\rangle=\langle 4,1\rangle=e,
\end{equation}
\begin{equation}\label{6.3}
[2,4]=[1,7]=e;~~1'=2~7~4~7~2;~~7'=7~4~1~2~1~4~7,
\end{equation}
\begin{equation}\label{6.4}
\langle 3,4\rangle=\langle 4,9\rangle=\langle 9,6\rangle=\langle 6,3\rangle=e,
\end{equation}
\begin{equation}\label{6.5}
[4,6]=[3,9]=e;~~3'=4~9~6~9~4;~~9'=9~6~3~4~3~6~9
\end{equation}
\begin{equation}\label{6.6}
\langle 5,6\rangle=\langle 6,10\rangle=\langle 10,8\rangle=\langle 8,5\rangle=e,
\end{equation}
\begin{equation}\label{6.7}
[6,8]=[5,10]=e;~~5'=6~10~8~10~6;~~10'=10~8~5~6~5~8~10,
\end{equation}
\begin{equation}\label{6.8}
1=3'~3~1'~3~3';~~~5~5'~5=3'~3~5~3~3';~~[1,5]=[1,5']=[1',5]=[1',5']=e,
\end{equation}
\begin{equation}\label{6.9}
\langle 7',9\rangle=\langle 7',9'\rangle=\langle 7',99'9\rangle=e;~~~\langle 9',10\rangle=\langle 10,9'99'7'9'99'\rangle=\langle 9,10\rangle=e,
\end{equation}
\begin{equation}\label{6.10}
7=9'~9~7'~9~9';~~10~10'~10=9'~9~10~9~9';~~[7,10]=[7,10']=[7',10]=[7',10']=e,
\end{equation}
\begin{equation}\label{6.11}
[i,j]=e, \ \  i=1,1', j=6,6',~8,8',~9,9',~10,10',
\end{equation}
\begin{equation}\label{6.12}
[i,j]=e, \ \  i=2,2', j=3,3',~5,5',~6,6',~8,8',~9,9',~10,10',
\end{equation}
\begin{equation}\label{6.13}
[i,j]=e, \ \  i=3,3', j=7,7',~8,8',~10,10';~~[i,j]=e, \ \  i=4,4', j=5,5'~8,8'~10,10',
\end{equation}
\begin{equation}\label{6.14}
[i,j]=e, \ \  i=5,5', j=7,7',~9,9'; ~~[i,j]=e, \ \  i=6,6', j=7,7',
\end{equation}
\begin{equation}\label{6.15}
[i,j]=e, \ \  i=7,7', j=8,8';~~[i,j]=e \ \  i=8,8', j=9,9',
\end{equation}
\begin{equation}\label{6.16}
10'\ 10\ 9'\ 9\ 8'\ 8\ 7'\ 7\ 6'\ 6\ 5'\ 5\ 4'\ 4\ 3'\ 3\ 2'\ 2\ 1'\ 1=e.
\end{equation}

Using the equalities from  (\ref{6.1}) and other simple relations, we obtain the following relations in the group $G_1$:
\begin{eqnarray}
&~j=j', \ \ j=1,2,3,\cdots,10, \\
&~2=1~7~4~7~1; \ 9=6~3~4~3~6; \ 10=8~5~6~5~8,\\
&~\langle 4,7\rangle=\langle 1,4\rangle=\langle 1,3\rangle=\langle 3,4\rangle=\langle 3,6\rangle=\langle 3,5\rangle=\langle 5,6\rangle=\langle 5,8\rangle=e,\\
&~[7,1]=[7,3]=[7,5]=[7,6]=[7,8]=[1,6]=[1,5]=[1,8]=e,\\
&~[4,6]=[4,5]=[4,8]=e,\\
&~[3,8]=[6,8]=e,\\
&~[1,4~3~4]=[6,3~5~3]=e.
\end{eqnarray}

\begin{figure}[ht]
\begin{center}
\scalebox{0.8}{\includegraphics{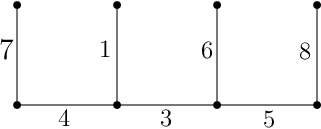}}
\end{center} \caption{The dual graph with generators related to $R_4\bigcup R_4$}\label{R4R4dual}
\end{figure}

From what has already been proved, we can choose $\{1,3,4,5,6,7,8\}$ as the generators of $G_1$. These choices give the graph of generators in Figure \ref{R4R4dual}. These relations show an isomorphism between $G_1$ and the symmetric group $S_8$. It follows that the Galois cover of surface $R_4 \bigcup R_4$ is a simply connected surface.
\end{proof}

\begin{figure}[ht]
\begin{center}
\scalebox{0.8}{\includegraphics{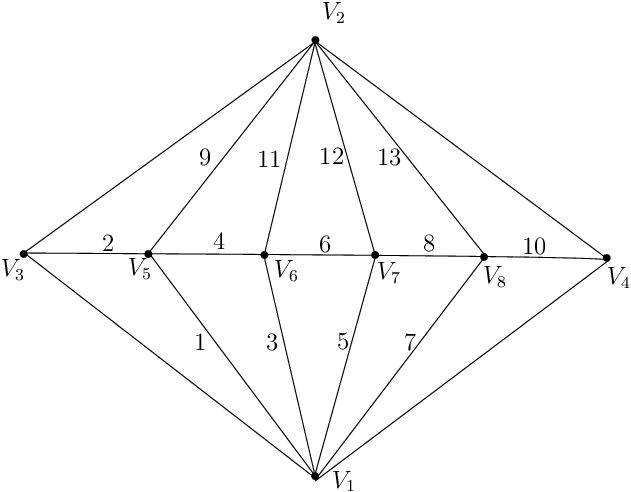}}
\end{center} \caption{The degeneration of the surface $R_5 \bigcup R_5$}\label{R5R5}
\end{figure}

We want to get a general result about the fundamental group of the Galois cover of a union of Zappatic surfaces, with each surface having one singularity of $R_{n+1}$. To obtain this final purpose, we must give a deeper discussion based on the above result and find some rules to solve it. In the following, it is natural to consider the case of a union of Zappatic surfaces, each of which with one singularity of type $R_{5}$, see Figure \ref{R5R5}. For this case, its result may be proved in much the same way as Theorem \ref{thm44}. Thus, we summarize without detailed proof the relevant material about the fundamental group of the Galois cover of a union of Zappatic surfaces, each of which with one singularity of type $R_{5}$. 

\begin{thm}\label{thm55}
The Galois cover of a union of Zappatic surfaces, each with one singularity of type $R_5$, is simply connected.
\end{thm}

\begin{proof}
In Figure \ref{R5R5}, it is clear that $V_1$ and $V_2$ are Zappatic singularities of type $R_5$, and their relations can be obtained from Lemma \ref{Rn-lemma1}. Each one of the vertices $V_i,\ i=5,6,7$, and $8$ is an intersection of four lines, and relations will come from Lemma \ref{basic-abcd}. The fundamental group $G=\pi_1(\mathbb{CP}^2-S_{R_5 \bigcup R_5})$ is generated by generators ${\{j,j'\}}^{13}_{j=1}$, but as for now, we work with $G_1=G/<j^2,j'^2>$.

Vertices $V_3$ and $V_4$ give rise to the relations $2=2'$, $10=10'$ in $G_1$. By Lemma \ref{innerR4-lemma1}, we get $i=i'$ for 
 $i=4,6,8$. By these results, we get a simplified presentation for $G_1$:

\begin{eqnarray}
&~j=j', \ j=1,2,3,\cdots,13, \\
&~2=1~9~4~9~1; \ 11=6~3~4~3~6; \ 12=8~5~6~5~8; \ 13=10~7~8~7~10,\\
&~\langle 4,9\rangle=\langle 1,4\rangle=\langle 1,3\rangle=\langle 3,4\rangle=\langle 3,6\rangle=e,\\
&\langle 3,5\rangle=\langle 5,6\rangle=\langle 5,8\rangle=\langle 5,7\rangle=\langle 8,7\rangle=\langle 7,10\rangle=e,\\
&~[9,1]=[9,3]=[9,5]=[9,6]=[9,7]=[9,8]=[9,10]=e,\\
&~[1,5]=[1,6]=[1,7]=[1,8]=[1,10]=[3,7]=[3,8]=[3,10]=e,\\
&~[4,5]=[4,6]=[4,7]=[4,8]=[4,10]=[6,7]=[6,8]=[6,10]=e,\\
&~[5,10]=[8,10]=e,\\
&~[1,4~3~4]=[6,3~5~3]=[8,5~7~5]=e.
\end{eqnarray}

\begin{figure}[ht]
\begin{center}
\scalebox{0.8}{\includegraphics{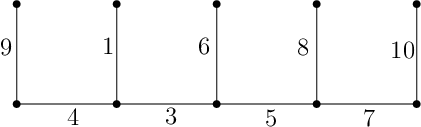}}
\end{center} \caption{The dual graph of the generators related to $R_5\bigcup R_5$}\label{R5R5dual}
\end{figure}

We choose $\{1,3,4,5,6,7,8,9,10\}$ as the generators of $G_1$. These choices provide the graph of generators in Figure \ref{R5R5dual} and conclude that $G_1\cong S_{10}$. Therefore, the Galois cover of surface ${R_5 \bigcup R_5}$ is a simply connected surface.
\end{proof}

\subsection{A union of Zappatic surfaces one singularity of type $R_{n+1}$}\label{n+1-type}

Having disposed of the above questions in Section \ref{6-type}, these results are likely to produce good results by simplifying relations and finding some rules. Now, we continue in this fashion, obtaining general result in this section. In the following, our main results are stated and proved.

Let $R_{n+1} \bigcup R_{n+1}$ be a union of two Zappatic surfaces, with each surface degenerating to one $R_{n+1}$ singularity. Let $S_{R_{n+1} \bigcup R_{n+1}}$ be the branch curve.

\begin{thm}\label{thmnn}
The Galois cover of a union of Zappatic surfaces, with each surface having one singularity of type $R_{n+1}$, is simply connected.
\end{thm}

\begin{proof}
In Figure \ref{RnRn}, $V_1$ and $V_2$ are Zappatic singularities, with each singularity being of type $R_{n+1}$ and each vertex $V_i, \ i=5,6,\dots,n+3,n+4$ is intersection of four lines. We have relations of these vertices. The fundamental group $G=\pi_1(\mathbb{CP}^2-S_{R_{n+1} \bigcup R_{n+1}})$ is generated by generators $\{j,j'\}^{3n+1}_{j=1}$.

Vertices $V_3$ and $V_4$ give rise to the relations $2=2'$, $2n+2=(2n+2)'$ in $G_1$.
Having relation $2=2'$ in $V_5$, we get relation $4=4'$ by Lemma \ref{innerR4-lemma1}. Similarly, we get from $V_i, i=4, 6, \dots , 2n+2$ that $i=i'$.

Lastly, we collect these resulting equalities, together with other relations, which means that relations related to $V_5, V_6, \dots, V_{n+4}$ will be listed here (using the conditions in Lemma \ref{innerR4-lemma1}), and relations related to $V_1, V_2$ will be listed here as well (Lemma \ref{Rn-lemma1}).

\begin{eqnarray}
& {\underline{Equalities}} \ \ \ 2=2'; \ 4=4';\dots;\ 2n=(2n)';\ 2n+2=(2n+2)'.\\
& {\underline{Vertex~V_5}} \ \ \ \langle 1,2\rangle=\langle 2,2n+1\rangle=\langle 2n+1,4\rangle=\langle 4,1\rangle=e,\\
&[2,4]=[1,2n+1]=e;~~1'=2~(2n+1)~4~(2n+1)~2,\\
&(2n+1)'=(2n+1)~4~1~2~1~4~(2n+1).
\end{eqnarray}

\begin{eqnarray}
& {\underline{Vertex~V_6}} \ \ \ \langle 3,4\rangle=\langle 4,2n+3\rangle=\langle 2n+3,6\rangle=\langle 6,3\rangle=e,\\
&[4,6]=[3,2n+3]=e;~~3'=4~(2n+3)~6~(2n+3)~4,\\
&(2n+3)'=(2n+3)~6~3~4~3~6~(2n+3).
\end{eqnarray}

\begin{equation*}
\cdots\cdots\cdots
\end{equation*}

\begin{eqnarray}
&{\underline{Vertex~V_{n+3}}} \ \ \ \langle 2n-3,2n-2\rangle=\langle 2n-2,3n\rangle=\langle 3n,2n\rangle=e,\\
&~\langle 2n,2n-3\rangle=e;~~[2n-2,2n]=[2n-3,3n]=e,\\
&~(2n-3)'=(2n-2)~(3n)~(2n)~(3n)~(2n-2),\\
&~(3n)'=(3n)~(2n)~(2n-3)~(2n-2)~(2n-3)~(2n)(3n).
\end{eqnarray}

\begin{eqnarray}
&{\underline{Vertex~V_{n+4}}} \ \ \ \langle 2n-1,2n\rangle=\langle 2n,3n+1\rangle=e,\\
&~\langle 3n+1,2n+2\rangle=\langle 2n+2,2n-1\rangle=e,\\
&~[2n,2n+2]=[2n-1,3n+1]=e,\\
&~(2n-1)'=(2n)~(3n+1)~(2n+2)~(3n+1)~(2n),\\
&~(3n+1)'=(3n+1)~(2n+2)~(2n-1)~(2n)~(2n-1)~(2n+2)~(3n+1).
\end{eqnarray}

\begin{eqnarray}
&{\underline{Vertex~V_1}} \ \ \ \langle 1',3\rangle=\langle 1',3'\rangle=\langle 1',33'3\rangle=e;~\langle 3',5\rangle=\langle 5,3'33'1'3'33'\rangle=\langle 3,5\rangle=e,\\
&\langle 5,7\rangle=\langle5',7\rangle=e; \cdots; \langle 2n-5,2n-3\rangle=\langle(2n-5)',2n-3\rangle=e,\\
&\langle 2n-3, 2n-1\rangle=\langle(2n-3)', 2n-1\rangle=e,\\
&~1=3'31'33';~~~55'5=3'3533';~~~77'7=5'57'55'; \cdots,\\
&(2n-3)~(2n-3)'~(2n-3)=(2n-5)'~(2n-5)~(2n-3)'~(2n-5)~(2n-5)',\\
& (2n-1)~(2n-1)'~(2n-1)=(2n-3)'~(2n-3)~n'~(2n-3)~(2n-3)',\\
&~[1,5]=[1,5']=[1',5]=[1',5']=e,\\
&~[1,7]=[1',7]=\cdots=[1,2n-3]=[1',2n-3]=[1, 2n-1]=[1', 2n-1]=e,\\
&[3,7]=[3',7]=\cdots=[3,2n-3]=[3',2n-3]=[3, 2n-1]=[3', 2n-1]=e,\\
&[5,9]=[5',9]=\cdots=[5,2n-3]=[5',2n-3]=[5, 2n-1]=[5', 2n-1]=e,\\
&\cdots\cdots\cdots\\
&[2n-5, 2n-1]=[(2n-5)', 2n-1]=e.
\end{eqnarray}

\begin{eqnarray}
&{\underline{Vertex~V_2}} \ \ \ \langle (2n+1)', 2n+3\rangle=\langle (2n+1)', (2n+3)'\rangle=e, \\
&\langle (2n+1)', (2n+3) (2n+3)'(2n+3)\rangle=e,\\
&~\langle (2n+3)', 2n+5\rangle=\langle 2n+5, (2n+3)'(2n+3) (2n+3)'(2n+1)',\\
&(2n+3)'(2n+3) (2n+3)'\rangle=e, \\
&\langle 2n+3, 2n+5\rangle=e,\\
&\langle 2n+5, 2n+7\rangle=\langle(2n+5)', 2n+7\rangle=\cdots; \langle 3n-1,3n\rangle=e, \\
&\langle(3n-1)',3n\rangle=\langle 3n, 3n+1 \rangle=\langle(3n)', 3n+1\rangle=e,\\
&~ 2n+1=(2n+3)'~(2n+3)~(2n+1)'~(2n+3)~(2n+3)';~~~\\
&(2n+5)~(2n+5)'~(2n+5)=(2n+3)'~(2n+3)~(2n+5)~(2n+3)~(2n+3)',\\
&(2n+7)~(2n+7)'~(2n+7)=(2n+5)'~(2n+5)~(2n+7)'~(2n+5)~(2n+5)'; \cdots,\\
&(3n)~(3n)'~(3n)=(3n-1)'~(3n-1)~(3n)'~(3n-1)~(3n-1)',\\
& (3n+1)~(3n+1)'~(3n+1)=(3n)'~(3n)~(3n+1)'~(3n)~(3n)',
\end{eqnarray}

\begin{eqnarray}
&~[ 2n+1, 2n+5]=[ 2n+1, (2n+5)']=[ (2n+1)', 2n+5]=[ (2n+1)', (2n+5)']=e,\\
&~[ (2n+1), 2n+7]=[ (2n+1)', 2n+7]=\cdots=[2n+1,3n]=[ (2n+1)',3n]=e,\\
&[2n+1, 3n+1]=[ (2n+1)', 3n+1]=e,\\
&[2n+3, 2n+7]=[ (2n+3)', 2n+7]=\cdots=[2n+3,3n]=[ (2n+3)',3n]=e,\\
&[2n+3, 3n+1]=[ (2n+3)', 3n+1]=e,\\
&[2n+5,2n+9]=[ (2n+5)',2n+9]=\cdots=[2n+5,3n]=[ (2n+5)',3n]=e,\\
& [2n+5, 3n+1]=[ (2n+5)', 3n+1]=e,\\
&\cdots\cdots\cdots\\
&[3n-1, 3n+1]=[(3n-1)', 3n+1]=e.
\end{eqnarray}

\begin{eqnarray}
& {\underline{Relations~of~type~(\ref{par})}} \ \ \ [i,j]=e,~~i=1,1', j=6,6'~8,8'~\cdots~(2n+2),(2n+2)',~\\
&(2n+3),(2n+3)',~\cdots~(3n+1),(3n+1)',\\
&~[i,j]=e,~~i=2,2', j=3,3',~5,5',~6,6',~7,7',~\cdots~,\\
&(2n),(2n)',~(2n+2),(2n+2)',~\cdots~(3n+1),(3n+1)',\\
&\cdots\cdots\cdots\\
&~[i,j]=e,~~i=(2n+1),(2n+1)', j=(2n+2),(2n+2)',\\
&~[i,j]=e,~~i=(2n+2),(2n+2)', j=(2n+3),(2n+3)',~\\
&(2n+4),(2n+4)',~\cdots~(3n),(3n)'.
\end{eqnarray}

\begin{equation}\label{eq3.10}
{\underline{Relations~of~type~(\ref{vK4})}} \ \ \ (3n+1)'~(3n+1)~(3n)'~(3n)~\cdots\cdots~2'~2~1'~1=e.
\end{equation}

Now, we simplify the above relations and obtain the following relations in the group $G_1$:

\begin{eqnarray}
&~j=j', \ j=1,2,3,\cdots,(3n+1), \\
&~2=1~(2n+1)~4~(2n+1)~1;~~(2n+3)=6~3~4~3~6;~~(2n+4)=8~5~6~5~8;\cdots,\\
&(3n+1)=(2n+2)~(2n-1)~(2n)~(2n-1)~(2n+2),\\
&~\langle 4,2n+1\rangle=\langle 1,4\rangle=\langle 1,3\rangle=\langle 3,4\rangle=\langle 3,6\rangle=\langle 3,5\rangle=\langle 5,6\rangle=\langle 5,8\rangle=\cdots=e, \\
&\langle 2n-3,2n-2\rangle=\langle 2n-3,2n-1\rangle=\langle 2n-3,2n\rangle=e, \\
&\langle 2n-1,2n\rangle=\langle 2n-1,2n+2\rangle=e,\\
&~[2n+1,1]=[2n+1,3]=[2n+1,5]=[2n+1,6]=e,\\
&[2n+1,7]=\cdots=[2n+1,2n]=[2n+1,2n+2]=e,\\
&~[1,5]=[1,6]=\cdots=[1,2n+2]=e,\\
&~[4,5]=[4,6]=\cdots=[4,2n+2]=e,\\
&~\cdots\cdots\cdots\\
&~[2n-3,2n+2]=[2n,2n+2]=e,\\
&~[1,4~3~4]=[6,3~5~3]=[8,5~7~5]=\cdots=e, \\
&[2n-2,(2n-5)~(2n-3)~(2n-5)]=[2n,(2n-3)~(2n-1)~(2n-3)]=e.
\end{eqnarray}

\begin{figure}[ht]
\begin{center}
\scalebox{0.9}{\includegraphics{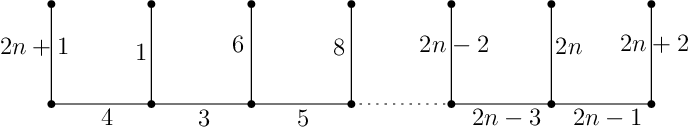}}
\end{center} \caption{The dual graph of $R_{n+1}\bigcup R_{n+1}$}\label{RnRndual}
\end{figure}

We choose $\{1,3,4,5,\cdots,2n+1,2n+2\}$ as the generators of $G_1$. This choice gives the dual graph of generators in Figure \ref{RnRndual}. The above-resulting relations can be mapped isomorphically onto the symmetric group $S_{2n+2}$; hence $G_1\cong S_{2n+2}$. It follows that the Galois cover of  $R_{n+1}\bigcup R_{n+1}$ is a simply connected surface.
\end{proof}

\section{Topological index}\label{Chern}

In this section, we will compute the Chern numbers of these Galois covers. The topological index $\tau(X_{Gal})=\frac{1}{3}(C_1^2(X_{Gal})-2C_2(X_{Gal}))$ of the Galois cover of a surface $X$ is an interesting invariant in the geography of algebraic surfaces (see \cite{Miy}). As a result, we will prove that these are negative indices. 
Let  $S\in \mathbb{CP}^2$ be the branch curve of $f$, where $f: X \rightarrow \mathbb{CP}^2$ is the generic projection of degree $N$.
In \cite{MoTe87}, we have the following formula:
$$ C_1^2(X_{Gal})=\frac{N!}{4}(m-6)^2$$
$$ C_2(X_{Gal})=N!(\frac{m^2}{2}-\frac{3m}{2}+3-\frac{3}{4}p-\frac{4}{3}q)$$
where $m$ is the degree of $S$, and $p$ is the number of nodes in $S$, $q$ is the number of cusps in $S$.

The values of $m$ and $N$ will be $6n+2$ and $2n$ respectively, according to Figure \ref{RnRn}. 
To understand the numerical values of $p,q$, we count all commutations (this will be $p$) and the number of relations of type $\langle * \ , \ * \rangle$ (this will be $q$). In \cite{AGM}, for a singularity of type $R_{n+1}$, we have $3(n-1)$ cusps and $(n-2) (2n-2)$ nodes, and from  Lemma \ref{basic-abcd} about a singularity with four intersecting lines, we have $12$ cusps and $4$ nodes.  

In this paper, for a surface $R_{n+1} \bigcup R_{n+1}$, we have 2 singularities of type $R_{n+1}$ and $n$ singularities that are intersections of four lines. In addition, we have relations of type (\ref{par}) that come from $14n^2-14n$ nodes. Therefore, we have a total number of nodes, which is $p=18n^2-22n+8$. In addition, we have $q=18n-6$  cusps.

Therefore, we have 
$$C_1^2(X_{Gal})=(2n+2)!~(9n^2-12n+4),~~ C_2(X_{Gal})=(2n+2)!~(\frac{9}{2}n^2-\frac{9}{2}n+4),~$$
and the Galois cover $X_{Gal}$ is a surface of general type by $C_1^2(X_{Gal})>0$.
We conclude that
$$\tau(X_{Gal})=\frac{1}{3}(C_1^2(X_{Gal})-2C_2(X_{Gal}))=\frac{1}{3}(2n+2)!~(-3n-4)<0.$$
This finding means that such Galois covers satisfy Bogomolov's conjecture.


\begin{thebibliography}{ABCD}


\bibitem{AGM} Amram, M., Gong, C., Mo, J-L., {\it On the Galois covers of degenerations of surfaces of minimal degree}, Math. Nachr., {\bf 296}(4), 2023, 1351--1365.
 
\bibitem{degree6} Amram, M., Gong, C., Sinichkin, U., Tan, S.-L., Xu, W.-Y., Yoshpe, M., {\it Fundamental groups of Galois covers of degree 6 surfaces}, Journal of Topology and Analysis, {\bf 15}, 2023,  593--613.



\bibitem{A-R-T}  Amram, M., Lehman, R., Shwartz, R., Teicher, M., {\it Classification of fundamental groups of Galois covers of surfaces of small degree degenerating to nice plane arrangements}, in Topology of algebraic varieties and singularities, {\bf 538}, Contemp. Math., 63--92, Amer. Math. Soc., Providence, RI, 2011.




\bibitem{Enri} Artal-Bartolo, E., {\it Topology of arrangements and position of singularities}, Annales de laFaculté des Sciences de Toulouse, {\bf 23}, 2014,  223--265.



\bibitem{C-C-F-M-5} Calabri, A., Ciliberto, C., Flamini, F., Miranda, R.,
{\it  On degenerations of surfaces}, arXiv:math/0310009v2, 2008.


\bibitem{C-C-F-M-4} Calabri, A., Ciliberto, C., Flamini, F., Miranda, R.,
{\it On the $K^2$ of degenerations of surfaces and the multiple point formula},
Ann. of Math., {\bf 165}(2), 2007, 335--395.



\bibitem{C1}
Catanese, F., {\it On the moduli spaces of surfaces of general type},  J. Diff. Geom., {\bf 19}, 1984, 483--515.

\bibitem{C2}
Catanese, F., {\it (Some) old and new results on algebraic surfaces}, First European Congress of Mathematics, Birkhauser Basel, 1994, 445--490.

\bibitem{CLM} Ciliberto, C., Lopez, A., Miranda, R., {\it Projective degenerations of $K3$ surfaces, Gaussian maps, and Fano threefolds}, Invent. Math., {\bf 114}, 1993, 641--667.

\bibitem{cog1} Cogolludo Agust{\'{\i}}n, J. I., {\it Braid Monodromy of Algebraic Curves}, Annales mathématiques Blaise Pascal, {\bf 18}, 2011, 141--209 . 


\bibitem{Gie}
Gieseker, D., {\it Global moduli for surfaces of general type}, Invent. Math., {\bf 41}, 1977, 233--282.



\bibitem{Hu1}
Gritsenko, V.A., Hulek, K., Sankaran, G.K., {\it The Kodaira dimension of the moduli of K3 surfaces},
Invent. Math., {\bf 169}, 2007, 519--567.


\bibitem{Hu2}
Hulek, K., Sankaran, G.K., {\it The Kodaira dimension of certain moduli spaces of abelian surfaces},
Compos. Math., {\bf 90}(1), 1994, 1--35.


\bibitem{KT}
Kulikov, V., Teicher, M., {\it 
Braid monodromy factorizations and diffeomorphism types}, Izv. Ross. Akad. Nauk Ser. Mat., {\bf 64}, 2000, 89–120.

\bibitem{Libgober1}
Libgober, A., {\it Fundamental groups of the complements to plane singular curves}, In Algebraic geometry,
Bowdoin, 1985 (Brunswick, Maine, 1985), pages 29--45. Amer. Math. Soc., Providence, RI, 1987.

\bibitem{Miy}
Miyaoka, Y.: {\it Algebraic surfaces with positive indices. In: Classification of Algebraic and Analytic Manifolds},Birkhäuser,  1982.


\bibitem{MoTe87}
Moishezon, B.,  Teicher, M., {\it Simply-connected algebraic surfaces of positive index}, Invent. Math.,
{\bf 89}, 1987, 601--643.




\bibitem{Paris} Paris, L., {\it On the fundamental group of the complement of a complex hyperplane arrangement}, 
Adv. Stud. Pure Math., {\bf 27}, 2000, 257--272.




\bibitem{vk}
van Kampen, E.R., {\it On the fundamental group of an algebraic curve}, Amer. J. Math., {\bf 55}, 1933, 255--260.

\bibitem{zg2}
Zappa, G., {\it Applicazione della teoria delle matrici di Veblen e di Poincar\'{e} allo studio delle superficie spezzate in sistemi di piani}, Acta Pont. Accad. Sci.,  {\bf 7}, 1943, 21--25.



\bibitem{Zappa} Zappa, G., {\it Su alcuni contributi alla conoscenza della struttura topologica delle superficie algebriche, dati dal metodo dello spezzamento in sistemi di piani}, Acta Pont. Accad. Sci., {\bf 7}, 1943, 4--8.







\end{thebibliography}
\end{document}